\documentclass[a4paper]{article}
\usepackage[margin=3cm]{geometry}
\usepackage{amsthm}
\usepackage{amsmath}
\usepackage{amssymb}
\usepackage[colorlinks,
urlcolor=blue!50!black,
citecolor=red!50!black,
linkcolor=blue!50!black]{hyperref}
\usepackage{graphicx}

\usepackage[UKenglish]{babel}
\usepackage[utf8]{inputenc}
\usepackage{float}
\usepackage{xcolor}
\usepackage{diagbox}
\usepackage{caption}
\usepackage{subcaption}
\usepackage{csquotes}
\usepackage{multicol}
\usepackage{multirow}
\usepackage{enumerate}
\usepackage{todonotes}

\theoremstyle{plain}
\newtheorem{theorem}{Theorem}

\newtheorem{proposition}[theorem]{Proposition}

\newtheorem{remark}[theorem]{Remark}

\theoremstyle{definition}
\newtheorem{definition}[theorem]{Definition}

\newcommand{\fd}{\mathfrak{d}}
\newcommand{\ve}{\varepsilon}
\newcommand{\R}{\mathbb{R}}
\newcommand{\N}{\mathbb{N}}

\newcommand{\p}{\mathbb{P}}

\title{Liouville quantum gravity metrics are not doubling}
\author{Liam Hughes}
\date{}
\begin{document}
	\maketitle
\begin{abstract}
\noindent We observe that non-doubling metric spaces can be characterized as those that contain arbitrarily large sets of approximately equidistant points and use this to show that, for $\gamma \in (0,2]$, the $\gamma$-Liouville quantum gravity metric is almost surely not doubling and thus cannot be quasisymmetrically embedded into any finite-dimensional Euclidean space. This generalizes the corresponding result of Troscheit \cite{troscheit} for the Brownian map (which is equivalent to the case $\gamma = \sqrt{8/3}$).
\end{abstract}
\section{Introduction}
For $\gamma \in (0,2]$, a \emph{$\gamma$-Liouville quantum gravity} (LQG) surface is a kind of random surface parametrized by a planar domain $D$ and formally described by the random metric tensor 
\begin{equation}\label{eq:lqg}
	e^{\gamma h(z)} \, (dx^2+dy^2)
\end{equation}
where $h$ is some form of the \emph{Gaussian free field} (GFF) on the domain $D$ and $dx^2+dy^2$ is the Euclidean metric. Such surfaces are known to describe scaling limits of certain random planar map models (see~\cite{ghs} and the references therein). 

Since $h$ is not sufficiently regular to be defined as a random function on $D$, but is instead a random distribution (in the sense of Schwartz), the formula (\ref{eq:lqg}) does not make literal sense, and thus LQG surfaces are constructed as limits of regularized versions of $e^{\gamma h(z)} dz$. For the \emph{subcritical} case $\gamma \in (0,2)$, such a construction given in~\cite{ds} yields the \emph{$\gamma$-LQG measure}~$\mu_h$, a random measure on~$D$, and the \emph{$\gamma$-LQG boundary length} $\nu_h$, a random measure on $\partial D$, both of which are encompassed by Kahane's~\cite{kahane} general theory of \emph{Gaussian multiplicative chaos}. (The regularization used for $\gamma \in (0,2)$ does not work for the \emph{critical} case $\gamma = 2$, for which the LQG measure was constructed in~\cite{lqg_crit}.)

LQG surfaces can also be defined as random metric spaces. It was proven in \cite{lqg1,lqg2,lqg3} for $\gamma = \sqrt{8/3}$, and subsequently in \cite{dddf,dfgps,gm} for $\gamma \in (0,2)$, that there exists a unique random metric $\fd_h$ associated to the GFF $h$ that satisfies a certain list of axioms associated with LQG (i.e., $\fd_h$ is a geodesic metric that transforms appropriately under affine coordinate changes and adding a continuous function to $h$, and is determined by~$h$ in a local way). The metric $\fd_h$ arises as a subsequential limit of \emph{Liouville first passage percolation}, a family of random metrics obtained from a mollified version of the GFF. Later, in~\cite{dg_supercrit} the LQG metric corresponding to the critical value $\gamma = 2$ was constructed and shown to be unique, as were \emph{supercritical} LQG metrics which correspond to complex values of $\gamma$ with $|\gamma|=2$.

Liouville quantum gravity surfaces are conformally covariant as metric measure spaces in the following sense: given a conformal map $\psi: \widetilde{D} \to D$, if we set
\begin{equation}\label{eq:reparam}
	Q = \frac{2}{\gamma} + \frac{\gamma}{2}, \quad \widetilde{h} = h \circ \psi + Q\log{|\psi'|},
\end{equation}
then by~\cite[Prop.~2.1]{ds}, we almost surely have $\mu_{\widetilde{h}} = \mu_h \circ \psi$ and (as long as $\psi$ extends by continuity to a homeomorphism between the closures of $\widetilde{D}$ and $D$ in the Riemann sphere) $\nu_{\widetilde{h}} = \nu_h \circ \psi$, and by~\cite[Thm~1.3]{gm_c}, we almost surely have $\fd_{\widetilde{h}} = \fd_h \circ \psi$.

We can then consider random metric measure spaces parametrized by $(D,h)$ with~$D$ a domain in $\mathbb{C}$ and $h$ some form of the GFF on $D$, with random quantum area and boundary length measures given by $\mu_h$ and $\nu_h$ and a random metric given by $\fd_h$, up to conformal reparametrization -- that is, we define a \emph{quantum surface} as an equivalence class of such pairs $(D,h)$ under the equivalence relation that identifies pairs related by conformal coordinate changes as described by (\ref{eq:reparam}). 

In~\cite{lqg2} it is shown that, in the particular case $\gamma = \sqrt{8/3}$, a certain kind of LQG surface called a \emph{quantum sphere} is almost surely isometric to another sphere-homeomorphic random object, the \emph{Brownian map} introduced by Le Gall~\cite{legall}, and further that this isometry almost surely pushes forward the LQG measure $\mu_h$ to the natural measure on the Brownian map, so that the quantum sphere and Brownian map are isomorphic as metric measure spaces. The law of the Brownian map is, intuitively, that of a ``uniform random element'' from the set of metric spaces that are homeomorphic to the sphere $S^2$, and it was proven independently by Le Gall~\cite{legall} and Miermont~\cite{miermont} that the Brownian map is the scaling limit of uniform random planar~quadrangulations.

The Brownian map can be constructed using a continuous process (the \emph{Brownian snake}) parametrized by the \emph{continuum random tree} (CRT), a random metric space constructed from the graph of a Brownian excursion by identifying points connected by horizontal line segments that stay underneath the graph. The CRT arises as the scaling limit of uniform discrete plane~trees. 

LQG is known to describe the scaling limit of certain discrete conformal embeddings of certain kinds of random planar map (e.g., the Cardy embedding of a uniform random triangulation~\cite{cardy_embedding} and the Tutte embedding of the mated-CRT map~\cite{tutte}). Given a quantum surface $\mathcal{S}$ and an embedding $\mathcal{S}\to D$ into the plane (i.e., a particular choice of parametrization $(D,h)$) obtained via such a scaling limit, one might therefore expect the embedding $\mathcal{S}\to D$ to retain some vestige of the conformality of the discrete embeddings. It is meaningless to ask directly whether the embedding $\mathcal{S}\to D$ is conformal, as the complex structure on~$\mathcal{S}$ comes from the embedding in the first place. However, since~$\mathcal{S}$ is a metric space, one could ask whether the embedding is \emph{quasisymmetric}. Quasisymmetric maps are embeddings of metric spaces in which the distortion of the metric is uniformly controlled; when both domain and target space are open subsets of~$\mathbb{R}^n$, locally quasisymmetric maps are equivalent to locally quasiconformal mappings.

In \cite{troscheit}, Troscheit proved that the continuum random tree and the Brownian map almost surely cannot be embedded quasisymmetrically into $\R^n$ for any $n$. The method was to show that those spaces have the property that for every $N$ one can find sets of~$N$ points all roughly equidistant from each other. This makes it impossible for such a space to be \emph{doubling}, i.e.\ for there to be a constant $M$ such that every open metric ball can be covered by at most $M$ open balls of half its radius. Moreover, whether a space is doubling (equivalently, whether its \emph{Assouad dimension} is finite) is preserved under quasisymmetric maps, so that no non-doubling space embeds quasisymmetrically, for instance, into $\mathbb{R}^n$. Since the equivalence with Brownian surfaces only holds for $\gamma = \sqrt{8/3}$, we will instead use GFF techniques to find approximately equidistant sets of points, generalizing the result of~\cite{troscheit} to all $\gamma \in (0,2]$:

\begin{theorem}\label{thm:main}
	Let $D \subseteq \mathbb{C}$ be a domain and $h$ some variant of the Gaussian free field on $D$. Let $\gamma \in (0,2]$ and let $\fd_h$ be the $\gamma$-LQG metric on $D$ associated to $h$. Then the metric space $(D,\fd_h)$ almost surely cannot be embedded quasisymmetrically into any doubling metric space (in particular, into $\mathbb{R}^n$, the sphere $S^n$, or any complete $n$-dimensional Riemannian manifold with non-negative Ricci curvature, for any $n\in\mathbb{N}$).
\end{theorem}
The proof is similar in spirit to Troscheit's: we show that a set of~$N$ roughly equidistant points exists at a given scale with fixed positive probability, then use the scaling properties of LQG to establish that there is almost surely \emph{some} scale at which such a set exists. The key ingredients are: a \emph{Girsanov-type result} for the GFF that implies that the metric behaves in a prescribed manner with positive probability; the \emph{locality} of the metric w.r.t.\ the GFF, which implies that one can detect such behaviour just by exploring the field locally; and the \emph{near-independence} of the GFF across scales, which along with locality allows us to upgrade the positive probability result to an almost sure one. We thus expect that the roadmap established in~\cite{troscheit} and here should allow one to obtain analogous results for other metric spaces constructed in a local manner from a sufficiently random object (Brownian motion in the cases of the CRT and Brownian map; the GFF here) satisfying some notion of near-independence in disjoint domains.
\subsection*{Acknowledgements}
This work was supported by the University of Cambridge Harding Distinguished Postgraduate Scholars Programme. I thank my PhD supervisor, Jason Miller, for suggesting the problem.
\section{Preliminaries}
\subsection{Quasisymmetric embeddings and Assouad dimension}
We first recall the definitions of quasisymmetric embedding and Assouad dimension.
\begin{definition}[Quasisymmetric embeddings]
	Let $(X,d_X)$ and $(Y,d_Y)$ be metric spaces and $f\colon X\to Y$ an injective function. Let $\Psi\colon (0,\infty)\to(0,\infty)$ be an increasing homeomorphism. Then $f$ is \textbf{$\Psi$-quasisymmetric} (equivalently, a \textbf{$\Psi$-quasisymmetric embedding}) if, for any three distinct points $x,y,z \in X$, we have
	\begin{equation}\label{eq:qs}
		\frac{d_Y(f(x),f(y))}{d_Y(f(x),f(z))} \le \Psi\left(\frac{d_X(x,y)}{d_X(x,z)} \right).
	\end{equation}
	We say $f$ is \textbf{quasisymmetric} (equivalently, a \textbf{quasisymmetric embedding}) if there exists some~$\Psi$ for which $f$ is $\Psi$-quasisymmetric.
\end{definition}
Recall that \emph{quasiconformal} maps between planar domains are intuitively those that send infinitesimal circles to infinitesimal ellipses of bounded eccentricity. We also define quasiconformality for embeddings between general metric spaces: with $X, Y, f$ as above, and $K\ge 1$, $f$ is \textbf{$K$-quasiconformal} if, for all $x\in X$, $$\limsup_{r\downarrow 0} \frac{ \sup\{ d_Y(f(x),f(y)) : d_X(x,y) \le r \} }{ \inf\{ d_Y(f(x),f(y)) : d_X(x,y) \ge r \}} \le K.$$

For open subsets of $\R^n$ with $n\ge 2$, locally quasisymmetric embeddings are equivalent to locally quasiconformal embeddings~\cite[Cor.~2.6]{vaisala_qs}. Indeed, the following holds for $n\ge 2$ and $D$ any domain in $\R^n$. Firstly~\cite[Thm~2.3]{vaisala_qs}, for each $\Psi$, if $f\colon D\to \R^n$ is a locally $\Psi$-quasisymmetric embedding, then $f$ is in fact $K$-quasiconformal for some $K\ge 1$ depending only on $\Psi$ and $n$. Conversely~\cite[Thm~2.4]{vaisala_qs}, for each $K\ge 1$, if $f\colon D\to \R^n$ is $K$-quasiconformal, and $x\in D$, $\alpha>1$, $r>0$ such that $B(x,\alpha r)\subseteq D$, then $f|_{B(x,r)}$ is $\Psi$-quasisymmetric for some $\Psi$ depending only on $K$, $n$ and $\alpha$. Similar results hold for smooth connected Riemannian manifolds~\cite[Thm~2.6]{afanaseva_bilet}.

We next introduce the Assouad dimension of a metric space. This is defined somewhat similarly to the upper box-counting dimension, but can be strictly greater -- intuitively, this happens when, in each covering by boxes of a given scale, disproportionately many boxes are required to cover certain particularly thick parts of the space.
\begin{definition}[Assouad dimension]\label{def:dima}
	Let $X$ be a metric space. For $E \subseteq X$, let $N_r(E)$ be the smallest possible cardinality of a set of open balls of radius $r$ that cover $E$. Then the \textbf{Assouad dimension} $\mathrm{dim}_\mathrm{A} \, X$ of $X$ is defined by
	\begin{equation*}\label{eq:dima}
		\mathrm{dim}_\mathrm{A} \, X := \inf\{ \alpha \ge 0: \exists C \in (0,\infty) \mathrm{ s.t. } \forall 0 < r < R, \forall x\in X, N_r(B(x,R)) \le C(R/r)^\alpha \}.
	\end{equation*}
\end{definition}
In \cite{troscheit} an alternative definition of $\mathrm{dim}_\mathrm{A}$ is used that only quantifies over $R < 1$. The Assouad dimension thus defined can be strictly smaller than the one defined in Def.~\ref{def:dima}, though they are equal when $X$ is compact, and both give dimension $n$ for $\mathbb{R}^n$. Our results and proofs apply regardless of which definition is used, but we use Def.~\ref{def:dima} since under this definition we have the equivalence (\cite[Thm~13.1.1]{fraser}) that $\mathrm{dim}_\mathrm{A} \, X < \infty$ if and only if~$X$ is a \textbf{doubling space}, i.e.\ there exists a finite constant $M$ such that any open ball in~$X$ can be covered by at most $M$ open balls of half its radius. (Under the other definition this equivalence fails; for example, it assigns dimension zero, rather than infinity, to the set of integer sequences in~$\ell^2$, which is not doubling.)

As observed by Coifman and Weiss~\cite[Ch.~III,~Lemma~1.1]{coifman_weiss}, a sufficient condition for a metric space $X$ to be doubling is the existence of a \emph{doubling measure}, that is a Borel measure $\mu$ on $X$ for which there is a constant $D>0$ such that, for all $x\in X$ and $r>0$,$$0<\mu(B(x,2r))\le D\mu(B(x,r))<\infty.$$A partial converse holds: whilst noting that $\mathbb{Q}$ is a doubling space for which there is no doubling measure, Assouad~\cite{assouad_pseudo} conjectured that every \emph{complete} doubling space has a doubling measure, which was proven by Luukkainen and Saksman~\cite{luukkainen} building on Vol'berg and Konyagin's~\cite{vk} proof for compact spaces. The Bishop--Gromov inequality~(\cite[\S11.10,~Corollary~3]{bishop_crittenden}; see also \cite[\S2.1]{gromov_betti}) straightforwardly implies that, for any complete Riemannian manifold with non-negative Ricci curvature, the measure given by the volume form is doubling, and thus such manifolds are doubling spaces. 
\subsection{Gaussian free field}
The \textbf{Gaussian free field} (GFF) is a two-time-dimensional analogue of Brownian motion. We first recall the definition of the \emph{zero-boundary} GFF from~\cite[Def.~2.10]{gffm} on a domain $D \subset \mathbb{C}$ with \emph{harmonically non-trivial} boundary (meaning that a Brownian motion started from $z\in D$ will almost surely hit $\partial D$). Let $H_s(D)$ be the set of smooth functions with compact support contained in $D$, equipped with the \emph{Dirichlet inner product}
\begin{equation*}
	(f,g)_\nabla = \frac{1}{2\pi} \int_D \nabla f(x) \cdot \nabla g(x) \, dx,
\end{equation*}
and complete this inner product space to a Hilbert space $H(D)$.
Taking an orthonormal basis $(\varphi_n)$ of $H(D)$ and letting $(\alpha_n)$ be i.i.d.\ $N(0,1)$ variables, we then define the \textbf{zero-boundary GFF in~$D$} as a random linear combination of elements of $H(D)$ given by
\begin{equation}\label{eq:gff}
	h = \sum_n \alpha_n \varphi_n.
\end{equation}
Though this sum does not converge pointwise or in the norm on $H(D)$, it can be shown (see~\cite[Prop.~2.7]{gffm}) that it \emph{does} converge almost surely in the space of distributions (so that the $L^2$ pairing $f\mapsto (h,f)$ is a continuous linear functional on $H_s(D)$) and in the fractional Sobolev space $H^{-\ve}(D)$ for each $\ve > 0$, and moreover that the law of the limit $h$ does not depend on the choice of basis~$(\varphi_n)$.

We can instead set $D$ to be all of $\mathbb{C}$; as in~\cite[\S2.2.1]{ig4}, we define the \textbf{whole-plane Gaussian free field} $h$ in the same way as above except that we consider $h$ \emph{modulo additive constant}. This means that we quotient out by an equivalence relation $\sim$ on the space of distributions defined such that $h_1 \sim h_2$ if and only if~$h_1-h_2$ is a constant distribution, meaning that there exists $a\in\R$ such that $(h_1,f)-(h_2,f) = a\int_{\mathbb{C}} f(z)\, dz$ for all $f\in H_s(\mathbb{C})$. An equivalent approach is to consider $(h,f)$ to be defined not for all~$f \in H_s(\mathbb{C})$ but only if $f \in H_{s,0}$, the subspace of~$H_s(\mathbb{C})$ consisting of those functions whose integral over $\mathbb{C}$ is zero. We say a random distribution $\widetilde{h}$ on $\mathbb{C}$ is a \textbf{whole-plane GFF plus a continuous function} if $\widetilde{h}$ can be coupled with a whole-plane GFF $h$ such that $\widetilde{h}-h$ is almost surely a continuous function on $\mathbb{C}$. From the expression (\ref{eq:gff}), it is straightforward to prove the Girsanov-type result that, if $f$ is a deterministic smooth compactly supported function, then the laws of $h$ and $h+f$ are mutually absolutely continuous (see~\cite[Prop.~3.4]{ig1}); we will make repeated use of this fact. Moreover, the laws of the restrictions of the whole-plane GFF and the zero-boundary GFF on a proper domain $D$ to a bounded subdomain $W \subset D$ with $\mathrm{dist}\, (W, \partial D) > 0$ are mutually absolutely continuous~\cite[Prop.~2.11]
{ig4}. Another useful property is that, since the Dirichlet inner product is conformally invariant, so is the GFF; in particular the GFF is invariant under translations and scalings (though only modulo additive constant in the whole-plane case).

As well as on $H_s(D)$, the pairing $(h,\cdot)$ can be defined for certain measures; in particular, for $D \neq \mathbb{C}$, the average $h_r(z)$ on the circle $\partial B(z,r)$ can be defined, whilst for $D = \mathbb{C}$ such circle averages are defined modulo additive constant, and indeed we can fix the additive constant by, for instance, requiring $h_1(0)=0$. Moreover, the process$$\{h_r(z)\colon z\in\mathbb{C}, r\in(0,\infty)\}$$has an almost surely continuous version~\cite[Prop.~3.1]{ds}, such that, for $z\in\mathbb{C}$ fixed, the process $(h_{e^{-t}}(z)-h_1(z))_{t\in\R}$ is a standard two-sided Brownian motion~\cite[Prop.~3.3]{ds}.
\subsection{Liouville quantum gravity metrics}
\cite[Thm~1.2]{gm} states that for $\gamma \in (0,2)$ there exists a measurable map $h \mapsto \fd_h$, from the space of distributions on $\mathbb{C}$ with its usual topology to the space of metrics on $\mathbb{C}$ that are continuous w.r.t.\ the Euclidean metric, that is characterized by satisfying the following axioms whenever $h$ is a whole-plane GFF plus a continuous function:
\begin{description}
	\item[Length space] Almost surely, the $\fd_h$-distance between any two points of $\mathbb{C}$ is given by the infimum of the $\fd_h$-lengths of continuous paths between the two points.
	\item[Locality] For a deterministic open set $U \subseteq \mathbb{C}$, define the \emph{internal metric} $\fd_h(\cdot,\cdot;U)$ of~$\fd_h$ on~$U$ by setting $\fd_h(z,w;U)$ to be the infimum of the $\fd_h$-lengths of continuous paths between $z$ and $w$ that are contained in $U$, for each pair of points $z,w\in U$. Then $\fd_h(z,w;U)$ is almost surely determined by $h|_U$.
	\item[Weyl scaling] Let $\xi = \gamma/d_\gamma$ where $d_\gamma$ is the fractal dimension defined in~\cite{dg}. Then for $f\colon\mathbb{C}\to\R$ continuous and $z,w\in\mathbb{C}$, define
	\begin{equation*}
		(e^{\xi f}\cdot \fd_h)(z,w) = \inf_P \int_0^{\mathrm{length}(P;\fd_h)} e^{\xi f(P(t))} \, dt
	\end{equation*}
	where the infimum is over continuous paths $P$ from $z$ to $w$ parametrized at unit $\fd_h$-speed. Then almost surely $e^{\xi f}\cdot \fd_h = \fd_{h+f}$ for all continuous $f$.
	\item[Affine coordinate change] For each fixed deterministic $r>0$ and $z\in\mathbb{C}$ it is almost surely true that, for all $u,v\in\mathbb{C}$,
	\begin{equation*}
		\fd_h(ru+z,rv+z) = \fd_{h(r\cdot+z)+Q\log{r}} (u,v),
	\end{equation*}
	where $Q=2/\gamma+\gamma/2$ is the constant appearing in the coordinate change rule (\ref{eq:reparam}).
\end{description}
The map $h\mapsto \fd_h$ is unique in the sense that for any two such objects $\fd$, $\widetilde{\fd}$, there is a deterministic constant~$C$ such that whenever $h$ is a whole-plane GFF plus a continuous function, almost surely we have $\fd_h = C\widetilde{\fd}_h$. This unique (modulo multiplicative constant) object is the \textbf{(whole-plane) $\gamma$-LQG metric}. Following~\cite{gm} we fix the constant so that the median distance between the left and right boundaries of $[0,1]^2$ is $1$ when $h$ is normalized so that $h_1(0)=0$.

Note that the scaling of $\fd_h$ is controlled by $\xi$, rather than $\gamma$. Indeed, since adding a constant~$C$ to~$h$ scales $\mu_h$ by $e^{\gamma C}$, it should be true that $\fd_h$ is scaled by~$e^{\xi C}$, where $\xi := \gamma/d_\gamma$ and $d_\gamma$ is the Hausdorff dimension of the $\gamma$-LQG metric. Though it would be circular to define $d_\gamma$ in this way, since it occurs in the Weyl scaling axiom, a definition of~$d_\gamma$ for each $\gamma \in (0,2)$ was obtained in~\cite{dg} by considering distances in certain discrete approximations of $\gamma$-LQG. A posteriori, it was shown in~\cite{gp1} that the Hausdorff dimension of the $\gamma$-LQG metric is indeed $d_\gamma$.

The existence and uniqueness of the $\gamma$-LQG metric were shown for the critical case $\gamma = 2$ in~\cite{dg_supercrit}, as well as the supercritical case corresponding to $\gamma \in \mathbb{C}$ with $|\gamma|=2$. We do not treat the supercritical case here, since such metrics have \emph{singular points} that are at distance $\infty$ from all other points and thus, while lower semicontinuous, fail to induce the Euclidean topology; however, the critical 2-LQG metric was shown to induce the Euclidean topology in~\cite{dg_crit} and satisfies the above axiomatic characterization just as in the subcritical case, except that $\xi$ needs to be replaced by $\xi_c := \lim_{\gamma \uparrow 2} \gamma/d_\gamma$.
\subsection{Notation}
For $z \in \mathbb{C}$ and $R_2 > R_1 > 0$, let $ \mathbb{A}_{R_1,R_2}(z)$ be the Euclidean annulus $B(z,R_2) \setminus \overline{B}(z,R_1)$.
\section{Non-doubling metric spaces}
We begin by giving an alternate characterization of non-doubling metric spaces (equivalently, those with infinite Assouad dimension) that we will verify for the LQG metric in order to rule out embeddability into $\R^n$. Namely, we observe that having infinite Assouad dimension is equivalent to containing arbitrarily large finite sets of points that are all approximately equidistant from each other, a characterization that does not seem to have appeared in previous literature.
\begin{definition} 
	Let $(X,d)$ be a metric space. Given $N \in \N$ and $K > 1$, we say that distinct points $x_1, \ldots, x_N \in X$ form an \textbf{$(N,K)$-clique} if 
	\begin{equation*}
		\max_{1\le i < j\le N} d(x_i,x_j) \le K \min_{1\le i < j\le N} d(x_i,x_j).
	\end{equation*}
	For $K > 1$ we say $(X,d)$ is \textbf{$K$-cliquey} if it contains an $(N,K)$-clique for each $N\in\N$.
\end{definition}
Instead of considering $(N,K)$-cliques, \cite{troscheit} considers ``approximate $N$-stars'' in which the $N$ points of a clique are also roughly equidistant from a central point that is closer to each outer point than the outer points are to each other. The proofs in both that paper and this one actually find approximate $N$-stars, but for our purposes the more simply defined $(N,K)$-cliques suffice, since quasisymmetric images of $K$-cliquey spaces must have infinite Assouad dimension:
\begin{proposition}\label{prop:cliquey}
	Let $(X,d_X)$ be a $K$-cliquey space for some $K>1$ and $f\colon(X,d_X)\to(Y,d_Y)$ a quasisymmetric mapping. Then $\mathrm{dim}_\mathrm{A} \, Y = \infty$.
\end{proposition}
\begin{proof}
	Suppose $(X,d_X)$ is $K$-cliquey and $f\colon(X,d_X)\to(Y,d_Y)$ is $\Psi$-quasisymmetric. Suppose also that $\mathrm{dim}_\mathrm{A}\, Y < \infty$, so that there exist $\alpha, C \in (0,\infty)$ for which $N_r(B(y,R)) \le C(R/r)^\alpha$ whenever $0<r<R$ and $y\in Y$. 
	
	Choose $N > 4^\alpha C(\Psi(K)^2+1)^\alpha$ and let $x_1, \ldots, x_N$ form an $(N,K)$-clique in $X$. Now by \eqref{eq:qs}, for $1 \le i,j,k \le N$ distinct we have
	\begin{equation}\label{eq:ijk}
		\frac{d_Y(f(x_i),f(x_j))}{d_Y(f(x_i),f(x_k))} \le \Psi\left( \frac{d_X(x_i,x_j)}{d_X(x_i,x_k)} \right) \le \Psi(K) \le \Psi(K)^2+1,
	\end{equation}
	since $x_1, \ldots, x_N$ form an $(N,K)$-clique and $\Psi$ is increasing. Applying the above twice, for $1 \le i,j,k,l \le N$ distinct we have
	\begin{equation}\label{eq:ijkl}
		\frac{d_Y(f(x_i),f(x_j))}{d_Y(f(x_k),f(x_l))} = \frac{d_Y(f(x_i),f(x_j))}{d_Y(f(x_i),f(x_k))} \cdot \frac{d_Y(f(x_i),f(x_k))}{d_Y(f(x_k),f(x_l))} \le \Psi(K)^2 \le \Psi(K)^2 + 1.
	\end{equation}
	\eqref{eq:ijk} and \eqref{eq:ijkl} together imply that $f(x_1), \ldots, f(x_N)$ form an $(N,\Psi(K)^2+1)$-clique. Now set $r = \frac{1}{2} \min_{1\le i < j \le N} d_Y(f(x_i),f(x_j))$ and $R = 2 \max_{1 \le i < j \le N}  d_Y(f(x_i),f(x_j))$. Then $B(f(x_1),R)$ contains all the $f(x_i)$ but no open ball of radius $r$ can contain more than one of the $f(x_i)$, so $N_r(B(f(x_1),R))\ge N> 4^\alpha C(\Psi(K)^2+1)^\alpha$. But $R/r \le 4(\Psi(K)^2+1)$ since the $f(x_i)$ form a $(\Psi(K)^2+1)$-clique, so this contradicts $N_r(B(y_1,R)) \le C(R/r)^\alpha$ and we have $\mathrm{dim}_\mathrm{A} \, Y = \infty$.
\end{proof}
In fact, being $K$-cliquey is equivalent to not being doubling (cf.~\cite[Prop.~2.7]{troscheit}):
\begin{proposition}\label{prop:equivalence}
	Let $(X,d)$ be a metric space. The following are equivalent:
	\begin{enumerate}[(i)]
		\item $\mathrm{\dim}_\mathrm{A} (X) = \infty$;
		\item $X$ is not a doubling space;
		\item $X$ is $K$-cliquey for some $K > 1$;
		\item $X$ is $K$-cliquey for every $K > 1$.
	\end{enumerate}
\end{proposition}
\begin{proof}
	\textbf{(i) $\Rightarrow$ (ii):} Contrapositively, if $X$ is a doubling space, then it is straightforward to show that $\mathrm{\dim}_\mathrm{A} (X)<\infty$ by iterating the operation of covering a ball with a fixed number of balls of half its radius (see~\cite[Thm~13.1.1]{fraser}).
	
	\textbf{(ii) $\Rightarrow$ (iii):} If $X$ is not doubling, then given any $N \in \mathbb{N}$ we can find $x\in X$ and $R>0$ such that $B(x,R)$ cannot be covered by less than $N$ balls of radius $R/2$. Let $x_1=x$ and construct $x_2, \ldots, x_N$ inductively so that $x_k \in B(x,R) \setminus \bigcup_{i=1}^{k-1} B(x_i,R/2)$ for $k=2, \ldots, N$ (possible by choice of $x$ and $R$). Now for $1\le i<j \le N$ we have$$R/2\le d(x_i,x_j) \le d(x_i,x)+d(x,x_j)< 2R,$$so the $x_i$ form an $(N,4)$-clique. Thus $X$ is 4-cliquey.
	
	\textbf{(iii) $\Rightarrow$ (iv):} If $K>1$ and $X$ is $K$-cliquey, then for any $N$ we can find an $(R(N),K)$-clique $x_1, x_2, \ldots, x_{R(N)}$ in $X$, where $R(N)$ is the $N$th Ramsey number. By definition of $R(N)$ such a clique contains $N$ points whose pairwise distances are either all in $[\min_{i<j} d(x_i,x_j), K^{1/2}\min_{i<j} d(x_i,x_j)]$ or all in $(K^{1/2}\min_{i<j} d(x_i,x_j),\max_{i<j} d(x_i,x_j)]$, in either case forming an $(N,K^{1/2})$-clique. Iterating this argument, we find that $X$ is $K^{1/4}$-cliquey, $K^{1/8}$-cliquey, and so on.
	
	\textbf{(iv) $\Rightarrow$ (i):} Apply Prop.~\ref{prop:cliquey} to the identity on $X$.
\end{proof}
\begin{remark}
	From Prop.\ \ref{prop:cliquey} and Prop.\ \ref{prop:equivalence} we deduce the well-known result that quasisymmetric images of non-doubling spaces are not doubling, and conversely (since the inverse of a $\Psi$-quasisymmetric bijection is $1/\Psi(1/\cdot)$-quasisymmetric) that quasisymmetric images of doubling spaces are doubling.
\end{remark}
We briefly observe another property that contrasts spaces of infinite and finite Assouad dimension. Given a metric space $(X,d_X)$ and $\beta \in (0,1)$ one can define the \textbf{$\beta$-snowflaking}~$d_X^\beta$ of~$d_X$ as the metric on $X$ given by $d_X^\beta(x,y) = d_X(x,y)^\beta$. 
The \emph{Assouad embedding theorem} \cite[Prop.~2.6]{assouad} states that for each $\alpha\in(0,\infty)$ and $\beta\in(0,1)$ there exists $n = n(\alpha,\beta) \in \N$ such that, if $(X,d)$ is a metric space such that $\mathrm{dim}_\mathrm{A}(X)=\alpha$, then there is a bi-Lipschitz embedding of $(X,d^\beta)$ into~$\R^n$. 
(Naor and Neiman~\cite{naor_neiman} later proved that one can choose $n = n(\alpha)$ such that $\R^n$ admits bi-Lipschitz embeddings of $(X,d^\beta)$ for all $\beta \in (1/2,1)$ and all $X$ with $\mathrm{dim}_\mathrm{A}(X)=\alpha$.) For spaces of infinite Assouad dimension, however, snowflaking does not facilitate bi-Lipschitz embeddings into~$\R^n$:
\begin{remark}
	Note that, for $\beta \in (0,1)$, if $(X,d_X)$ is $K$-cliquey then $(X,d_X^\beta)$ is $K^\beta$-cliquey; thus, if $\mathrm{dim}_\mathrm{A} X = \infty$ then the $\beta$-snowflaking of $X$ cannot be embedded quasisymmetrically into any doubling space (cf.~\cite[Thm~2.6]{troscheit}).
\end{remark}
\section{Proof of Theorem \ref{thm:main}}
We will begin by proving the result for the whole-plane GFF and deduce it for other variants via local absolute continuity. The main task is to show that for $N$, $\delta$ fixed, a fixed closed disc contains an $(N,1+\delta)$-clique with positive probability (by scale and translation invariance, this probability will not depend on the disc). The basic idea for this is to consider a polygonal star with $N$ arms and add bump functions to the field in order to force geodesics between the arms to stay within the star, recalling that the law of the modified field will be mutually absolutely continuous with that of the original field. The near-independence of the field in disjoint regions then allows us to translate positive probability for a fixed disc into an almost sure result: a Markovian exploration of the domain (we will use the annulus exploration from \cite[Lemma~3.1]{gm_l}) will almost surely find a disc containing an $(N,1+\delta)$-clique. Since this holds for every $N$, we have that $\gamma$-LQG metric spaces are $(1+\delta)$-cliquey, so by Prop.~\ref{prop:cliquey} their quasisymmetric images have infinite Assouad dimension, which as mentioned is equivalent to not being doubling.

Fix $N \ge 2$, $z_0 \in\mathbb{C}$, $r > 0$, $\delta \in (0,1)$ and $\ve\in (0,1/14)$. Let $h$ be a whole-plane GFF, normalized so that (say) $h$ has average zero on some fixed circle disjoint from $\overline{B}(z_0,8r)$. Now set $z_k = z_0 + 6re^{2\pi ik/N}$, $z'_k = z_0 + 7re^{2\pi ik/N}$ and $w_k = z_0 + r e^{\pi i(2k+1)/N}$ for $k = 1, \ldots, N$, and let $K_N$ be the compact set consisting of the polygon whose sides are the line segments joining$$(z'_1, w_1), (w_1, z'_2), (z'_2, w_2), (w_2, z'_3), \ldots, (z'_N, w_N), (w_N, z'_1)$$together with this polygon's interior. For $\beta \in (0,1)$ let $K_N^\beta = z_0 + (1-\beta)(K_N-z_0)$.
\begin{figure}
	\centering
	\includegraphics[width=0.66\linewidth]{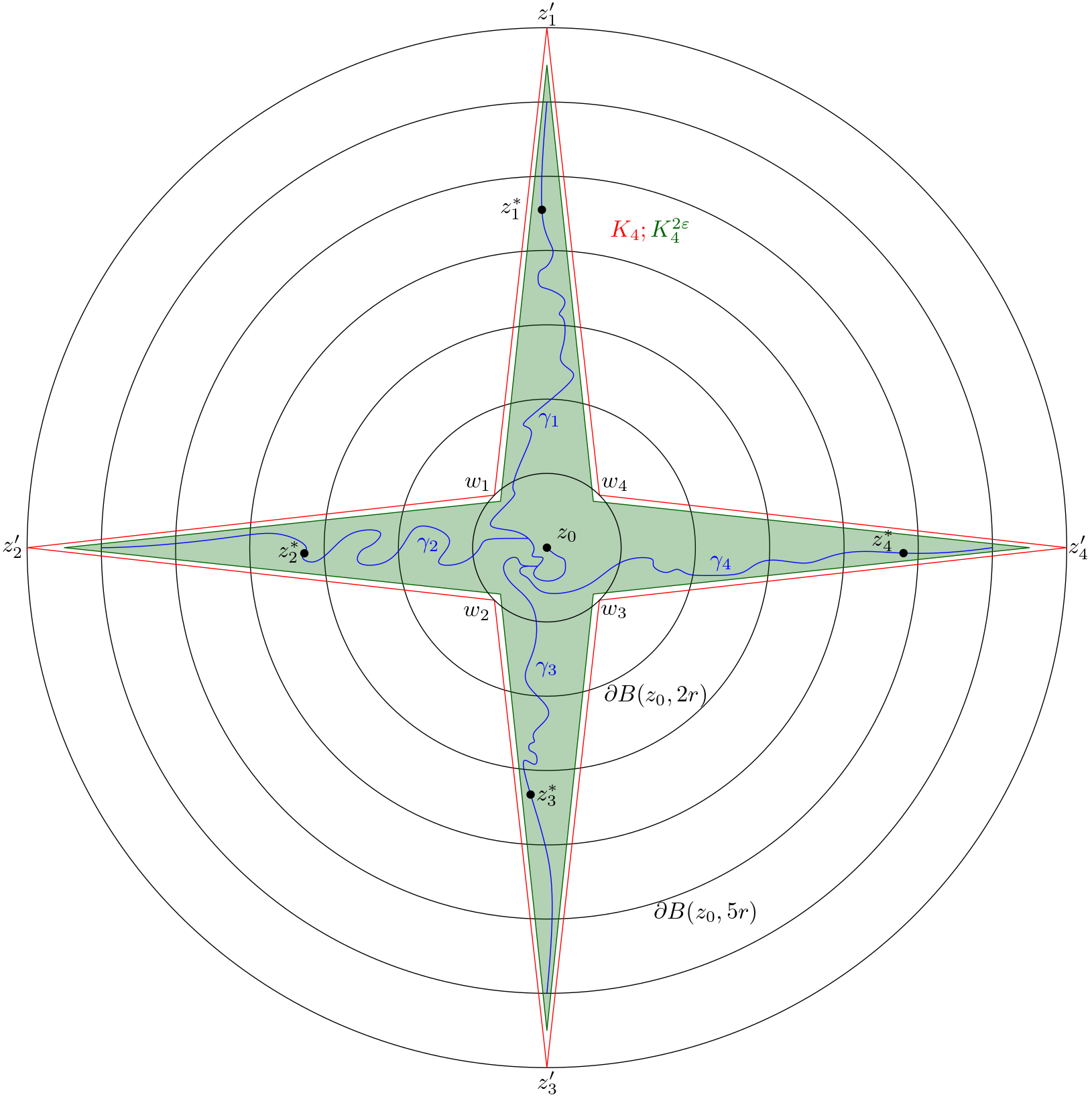}
	\caption{\label{fig:star}The points $z^*_i$ are chosen to be equidistant from $\partial B(z_0,2r)$; we then arrange that geodesics between them stay within $K_N \cup \overline{B}(z_0,2r)$ and that the diameter of $\overline{B}(z_0,2r)$ is small, making the $z^*_i$ almost equidistant from each other.}
\end{figure}
Fix $\zeta(\ve)>0$ such that the Euclidean $2\zeta(\ve)$-neighbourhood of $K^{\ve/2}_N$ is contained in $K^{\ve/4}_N$. Define the event$$A^1_C(h) = \begin{Bmatrix}
	\inf\left\{ \fd_h(z,w) \, \middle|  \, z,w \in B(z_0,7r)\setminus K^{\ve/2}_N, |z-w|\ge \zeta(\ve) \right\} > 1/C; \\
	\fd_h(\partial B(z_0,2r), \partial B(z_0,5r);\,\mathrm{int}\, K^{\ve}_N) \le C
\end{Bmatrix}. $$
If $\widetilde{h}$ is another field, we define $A^1_C(\widetilde{h})$ to be the event given by replacing $h$ by $\widetilde{h}$ throughout in the definition of $A^1_C(h)$. (In the rest of this proof, we will tacitly use further definitions of this kind.)

We check that $\p[A^1_C(h)]\to 1$ as $C\to\infty$. It suffices to observe that, almost surely,$$0<\fd_h(\partial B(z_0,2r), \partial B(z_0,5r))<\fd_h(\partial B(z_0,2r), \partial B(z_0,5r);\,\mathrm{int}\, K^{\ve}_N)<\infty, $$whilst$$\inf\left\{ \fd_h(z,w) \, \middle| \, z,w \in B(z_0,7r)\setminus K^{\ve/2}_N, |z-w|\ge \zeta(\ve) \right\}>0,$$ since if not we could find $z_{(n)},w_{(n)} \in B(z_0,7r)\setminus K^{\ve/2}_N$ for each $n\in\N$, with $|z_{(n)}-w_{(n)}|\ge \zeta(\ve)$ and $\fd_h(z_{(n)},w_{(n)})\to 0$ as $n\to\infty$, and by Bolzano--Weierstrass and continuity of $\fd_h$ w.r.t.\ the Euclidean metric, extract subsequences converging to $z$ and $w$ (w.r.t.\ both the Euclidean metric and $\fd_h$) with $\fd_h(z,w)=0$ but $|z-w|\ge \zeta(\ve)$, a contradiction. (We could also use the local H\"{o}lder continuity of the Euclidean metric w.r.t.\ $\fd_h$ as proven for the subcritical case in~\cite[Prop.~3.18]{dfgps} and extended to the critical case in~\cite[Prop.~1.10]{pfeffer2024weak}.) 

Since $\p[A^1_C(h)]\to 1$ as $C\to\infty$, we can choose $C_1>0$ such that $\p[A^1_{C_1}(h)]>0$. Let~$\psi$ be a bump function supported in $B(z_0,8r)\setminus K^{\ve}_N$ such that $\psi\equiv 1$ on $B(z_0,7r)\setminus K^{\ve/2}_N$. For $\eta>0$, let $E_\eta(h)$ be the event that
\begin{align*}\inf &\left\{ \fd_h(z,w;\mathbb{A}_{2r,7r}(z_0)\setminus K^{\ve/2}_N)\, \middle| \, z,w \in \mathbb{A}_{2r,7r}(z_0)\setminus K^{\ve/2}_N, |z-w|\ge \zeta(\ve) \right\} \\
	&\ge 2\left(\fd_h(\partial B(z_0,2r), \partial B(z_0,5r))+\eta\right).\end{align*}
If we choose $M$ such that $e^{\xi M} > 2C_1^2$, and choose $\eta < e^{\xi M}/(2C_1)-C_1$, then by Weyl scaling we have $A^1_{C_1}(h) \subseteq E_\eta(h+M\psi)$. Thus, since $h$ and $h+M\psi$ have mutually absolutely continuous laws and $\p[A^1_{C_1}(h)]>0$, we can fix $\eta_1>0$ so that $\p[E_{\eta_1}(h)]>0$.

Since $\fd_h(\partial B(z_0,2r),\partial B(z_0,5r))>0$ almost surely, we can fix $t=t(\eta_1)>0$ so that$$\p[\fd_h(\partial B(z_0,2r),\partial B(z_0,5r)) > t/(2\delta)]>1-\p[E_{\eta_1}(h)].$$Next, note that $\sup_{z\in \mathbb{A}_{(2-u)r,2r}(z_0)} \fd_h(z,\partial B (z_0,(2-u)r)) \to 0$ almost surely as $u\downarrow 0$, since if not we could find $v>0$ and sequences $z_{(n)}$, $u_{(n)}$ such that $u_{(n)} \downarrow 0$, $z_{(n)} \in \mathbb{A}_{(2-u_{(n)})r,2r}(z_0)$ and $\fd_h(z_{(n)},\partial B (z_0,(2-u_{(n)})r))\ge v$; then any subsequential limit $z\in\partial B(z_0,2r)$ must have $\fd_h$-distance $\ge v$ from $B(z_0,2r)$, a contradiction (indeed, given any $u \in (0,2)$, once~$n$ is large enough that $u_{(n)} \le u$ we have $\fd_h(z_{(n)},\overline{B}(z_0,(2-u)r)) \ge v$ and, taking the subsequential limit, $\fd_h(z,\overline{B}(z_0,(2-u)r)) \ge v$). This convergence holds almost surely and hence in probability, so we can fix $u=u(t,\eta_1)>0$ such that
\begin{align*}
	&\p \left[ \sup_{z\in\mathbb{A}_{(2-u)r,2r}(z_0)} \fd_h (z,\partial B (z_0,(2-u)r)) < t/3 \right] \\
	&> 1-\p\big[ \{ \fd_h(\partial B(z_0,2r),\partial B(z_0,5r)) > t/(2\delta) \} \cap E_{\eta_1} (h)\big]
\end{align*}
(since the right-hand side is strictly less than 1 by choice of $t$).  Since as $C\to\infty$ we have $\p[\mathrm{diam}\,(\overline{B}(z_0,(2-u)r);\fd_h(\cdot,\cdot;B(z_0,(2-u/2)r))) \le C]\to 1$, we can fix $C_2=C_2(u,t,\eta_1)>0$ so that
\begin{align*}
	\p[F_{C_2,t,u}(h)] &:=\p\begin{bmatrix}
	\sup_{z\in\mathbb{A}_{(2-u)r,2r}(z_0)} \fd_h(z,\partial B (z_0,(2-u)r)) < t/3; \\
	\mathrm{diam}\,(\overline{B}(z_0,(2-u)r);\fd_h(\cdot,\cdot;B(z_0,(2-u/2)r))) \le C_2
\end{bmatrix}\\
 &> 1-\p\big[\{ \fd_h(\partial B(z_0,2r),\partial B(z_0,5r)) > t/(2\delta) \} \cap E_{\eta_1(h)}\big],
\end{align*}i.e., with positive probability both $E_{\eta_1}(h)\cap F_{C_2,t,u}(h)$ holds and $\fd_h(\partial B(z_0,2r),\partial B(z_0,5r))>t/(2\delta)$. Now fix a bump function $\sigma$ such that $\sigma\equiv 1$ on $B(z_0,(2-u/2)r)$ and $\sigma\equiv 0$ outside $B(z_0,2r)$. For $M'$ large enough depending on $C_2$, $t$, $u$, on $F_{C_2,t,u}$ we have$$ \mathrm{diam}\,(\overline{B}(z_0,(2-u)r);\fd_{h-M'\sigma}) \le t/3, \quad \sup_{z\in\mathbb{A}_{(2-u)r,2r}(z_0)} \fd_{h-M'\sigma}(z,\partial B (z_0,(2-u)r)) < t/3$$(the latter because $\fd_{h-M'\sigma}\le\fd_{h}$ pointwise), which implies $\mathrm{diam}\, (\overline{B}(z_0,2r);\fd_{h-M'\sigma})<t$. Thus, on the event that both $E_{\eta_1}(h)\cap F_{C_2,t,u}(h)$ holds and $\fd_h(\partial B(z_0,2r),\partial B(z_0,5r))>t/(2\delta)$, we have $\mathrm{diam}\, (\overline{B}(z_0,2r);\fd_{h-M'\sigma})<t$, whilst $\fd_{h-M'\sigma}(\partial B(z_0,2r),\partial B(z_0,5r))>t/(2\delta)$, and the event $E_{\eta_1}(h-M'\sigma)$ holds. Indeed, the latter two events only depend on the field outside $B(z_0,2r)$ so are invariant under replacing $h$ by $h-M'\sigma$. Since $h$ and $h-M'\sigma$ have mutually absolutely continuous laws, we may conclude that with positive probability, $\mathrm{diam}\, (\overline{B}(z_0,2r);\fd_h)<t$, $E_{\eta_1}(h)$ holds, and $\fd_h(\partial B(z_0,2r),\partial B(z_0,5r))>t/(2\delta)$.

Since $\ve < 1/14$, we have $z_i \in \, \mathrm{int}\, K^{2\ve}_N$ for each $i$, $0\le i\le N$, so we can almost surely find paths $\gamma_i = \gamma_i(h|_{B(z_0,8r)}) \subset \mathrm{int}\, K^{2\ve}_N$ from $z_0$ to $z_i$ for $1 \le i \le N$ with finite $\fd_h$-length (e.g., by~\cite[Prop.~3.9]{dfgps}), which we can fix in some manner that is measurable w.r.t.~$h|_{B(z_0,8r)}$ considered modulo additive constant. For concreteness, we choose $m\in\mathbb{N}$ to be the (random) least number such that, whenever $1\le i\le N$, $1\le k\le m$, $\fd_h(\frac{1}{m}((m-k)z_0+kz_i),\frac{1}{m}((m-k+1)z_0+(k-1)z_i)) < \fd_h(\frac{1}{m}((m-k)z_0+kz_i),\partial K^{2\ve}_N)$, then choose $\gamma_i$ to be the concatenation of the almost surely unique $\fd_h$-geodesics from $\frac{1}{m}((m-k+1)z_0+(k-1)z_i)$ to $\frac{1}{m}((m-k)z_0+kz_i)$.

For $1\le i \le N$, explore $\gamma_i$ from $z_0$ towards $z_i$. By continuity of $\fd_h(\cdot,\partial B(z_0,2r))$, we can define~$z^*_i$ as the first point of $\gamma_i\setminus \overline{B}(z_0,2r)$ reached by this exploration such that$$\fd_h(z^*_i,\partial B(z_0,2r)) = \fd_h(\partial B(z_0,2r),\partial B(z_0,5r)).$$We argue that, on an event which we have just shown to have positive probability, namely$$G_{z_0,8r}(h) := E_{\eta_1}(h) \cap \{ \mathrm{diam}\, (\overline{B}(z_0,2r);\fd_h)<2\delta\fd_h(\partial B(z_0,2r),\partial B(z_0,5r)) \},$$the $z^*_i$ form an $(N,1+\delta)$-clique. On $E_{\eta_1}(h)$, for $1\le i<j \le N$ we have$$\fd_h(z^*_i,z^*_j) \ge \fd_h(z^*_i,\partial B(z_0,2r)) + \fd_h(z^*_j,\partial B(z_0,2r)) = 2\fd_h(\partial B(z_0,2r),\partial B(z_0,5r)).$$Indeed, this lower bound certainly holds for any path from $z^*_i$ to $z^*_j$ that intersects $\overline{B}(z_0,2r)$. Also, since $\overline{B}(z_0,2r)$ disconnects the prongs of the star $K_N^{\ve/4}$, any path~from~$z^*_i$ to~$z^*_j$ that does \emph{not} enter $\overline{B}(z_0,2r)$ must have a subpath contained in $\mathbb{A}_{2r,7r}(z_0)\setminus K_N^{\ve/2}$ of Euclidean diameter at least $\zeta(\ve)$, which on $E_{\eta_1}(h)$ must have $\fd_h$-length greater than $2\fd_h(\partial B(z_0,2r),\partial B(z_0,5r))$. 

Finally, on the event that $\mathrm{diam}\, (\overline{B}(z_0,2r);\fd_h)<2\delta\fd_h(\partial B(z_0,2r),\partial B(z_0,5r))$, we have
\begin{align*}\fd_h(z^*_i,z^*_j) &\le \fd_h(z^*_i,\partial B(z_0,2r)) + \fd_h(z^*_j,\partial B(z_0,2r)) + \mathrm{diam}\, (\overline{B}(z_0,2r);\fd_h) \\
	&< 2(1+\delta)\fd_h(\partial B(z_0,2r),\partial B(z_0,5r)).
\end{align*}
Therefore the $z^*_i$ form an $(N,1+\delta)$-clique. Thus, on $G_{z_0,8r}(h)$, there exist points in $B(z_0,8r)$ that form an $(N,1+\delta)$-clique w.r.t.\ $\fd_h$. 

Note that, since $G_{z_0,8r}(h)$ only depends on ratios between distances and thus is determined by the field modulo additive constant, the scale and translation invariance properties of $h$ imply that the analogous event $G_{z,r'}(h)$ with $z_0$ and~$r$ replaced respectively by $z$ and $r'$ (and the necessary changes made in the definitions of $\gamma_i$, $z^*_i$, $E_{\eta_1}(h)$) has the same probability for any $z \in \mathbb{C}$ and any $r' > 0$. Moreover, since $G_{z,r'}(h)$ is determined by $h|_{B(z,8r')}$, it is in fact determined by $(h-h_R(z'))|_{B(z,8r')}$ whenever $\overline{B}(z,8r')\subset B(z',R)$.

We can now consider a sequence of nested concentric annuli within which we have near-independence of the field, meaning that if we take a closed disc $\overline{B} (z^{(k)},8r^{(k)})$ within each annulus then at least one of the events $G_{z^{(k)},r^{(k)}}(h)$ holds. Indeed we are in the setting of~\cite[Lemma~3.1]{gm_l}, which implies that, say, for the annuli $(\mathbb{A}_{2^{-2k-1},2^{-2k}}(0))_{k\in \mathbb{N}}$ and $z^{(k)} = 3\cdot 2^{-2k-2}$, $8r^{(k)}=2^{-2k-3}$, a positive proportion of the events $\{ G_{z^{(k)},8r^{(k)}}(h) \}_{k=1}^K$ hold with probability exponentially high in $K$. In particular, it is almost surely the case that at least one of the events $G_{z^{(k)},8r^{(k)}}(h)$ holds.

Since we have now shown that an $(N,1+\delta)$-clique almost surely exists for all $N$ within a fixed closed disc, the surface $(\mathbb{C},\fd_h)$ is almost surely $(1+\delta)$-cliquey and thus cannot be embedded quasisymmetrically into any doubling space. The fact that this argument finds all the $(N,1+\delta)$-cliques within the same disc also means that the local mutual absolute continuity of GFF variants gives the same result for other LQG surfaces, and thus we conclude the proof of Theorem~\ref{thm:main}.
\bibliographystyle{alpha}
\bibliography{newbib}
\end{document}